\documentclass[11pt]{amsart}
\usepackage{graphicx,amssymb,amsmath,epsfig}
\newtheorem{theorem}{Theorem}[section]
\newtheorem{proposition}[theorem]{Proposition}
\newtheorem{lemma}[theorem]{Lemma}
\newtheorem{corollary}[theorem]{Corollary}
\newtheorem{remark}[theorem]{Remark}
\newtheorem{example}[theorem]{Example}

\newtheorem{definition}[theorem]{Definition}
\newtheorem{properties}[theorem]{Properties}\rm

\usepackage{amssymb,amsmath}
\usepackage[all]{xy}
\usepackage{enumerate}

\newcommand{\bth}{\begin{theorem}}
\newcommand{\bpr}{\begin{proposition}}
\newcommand{\epr}{\end{proposition}}
\newcommand{\bco}{\begin{corollary}}
\newcommand{\eco}{\end{corollary}}
\newcommand{\ble}{\begin{lemma}}
\newcommand{\ele}{\end{lemma}}
\newcommand{\bre}{\begin{remark}\rm}
\newcommand{\ere}{\end{remark}}
\newcommand{\bex}{\begin{example}\rm}
\newcommand{\eex}{\end{example}}
\newcommand{\bde}{\begin{definition}\rm}
\newcommand{\ede}{\end{definition}}

\def\la#1{\hbox to #1pc{\leftarrowfill}}
\def\ra#1{\hbox to #1pc{\rightarrowfill}}
\def\fract#1#2{\raise4pt\hbox{$ #1 \atop #2 $}}

\def\lrar{{\longrightarrow}\ }

\def\bbz{{\mathbb Z}}

\def\bbr{{\mathbb R}}

\def\bbc{{\mathbb C}}
\def\bbr{{\mathbb R}}

\def\map#1{\hbox{Map}_{#1}}

\begin{document}

\title{ Thin Loop Groups}%
\author{Moncef Ghazel}\author{Sadok Kallel}
\address{ Facult\'e des Sciences Math\'ematiques, Physiques et Naturelles de Tunis. Universit\'e de Tunis El Manar.}
\address{American University of Sharjah (UAE) and Laboratoire Painlev\'e, Lille 1 (France).}
\email{sadok.kallel@univ-lille.fr}

\maketitle

\centerline{In fond memory of Stefan Papadima}

\begin{abstract}
We verify that for a finite  simplicial complex $X$ and for piecewise linear loops on $X$,  the ``thin" loop space is a topological group of the same homotopy type as the space of continuous loops. This turns out not to be the case for the higher loops.
\end{abstract}

\section{Introduction}\label{intro}

A loop on a topological space $X$ is any continuous map of the circle $S^1\subset\bbc$ into $X$. A loop is ``based" if it sends $1\in S^1$ to a preferred basepoint $x_0\in X$.
Physicists and geometers as in \cite{barrett, caetano} introduce the ``thin
fundamental group'' of a smooth manifold to be the space of based smooth loops, ``stationary" at the basepoint, modulo ``thin homotopy''. A homotopy is ``thin'' if
its image ``doesn't enclose area'' (precise definitions in \S\ref{thin}).
If $X$ is a smooth manifold, then its thin fundamental group is denoted by $\pi_1^1(X)$. This is indeed a group which is a quotient of the corresponding loop space $\Omega (X)$ and of which $\pi_1(X)$ is a quotient; i.e.
\begin{equation}\label{sequence}
\Omega (X)\twoheadrightarrow\pi_1^1(X)\twoheadrightarrow\pi_1(X)
\end{equation}
This space plays an important role in studying holonomy, connections and transport.
As shown in \cite{barrett, picken}, and even earlier in \cite{kobayashi, teleman}, there is a bijective correspondence between homomorphisms
$\xymatrix{\pi_1^1(X)\ar[r]^h&G}$ with $G$ a compact Lie group, up to conjugation, and principal $G$-bundles with connections on $X$, unique up to equivalence. A technical condition is put on the topology of loops so that there is the notion of a smooth family of loops which has a smoothly varying holonomy image in $G$ (\cite{barrett},\S2). The interesting part of this observation is that if $f,g\in \Omega (X)$ are two smooth loops related by a thin homotopy, then parallel transport along $f$ gives the same result as parallel transport along $g$ \cite{caetano, kobayashi}.

The question we ask, and then answer in this paper, is: how ``close" is $\pi_1^1(X)$ to $\pi_1(X)$ or to $\Omega (X)$, the space of all continuous loops, when $X$ is a finite complex? When $X$ is a graph or a manifold of dimension $1$, all homotopies are thin so the epimorphism
$\pi_1^1(X)\twoheadrightarrow\pi_1(X)$
becomes an isomorphism. This situation turns out to be rather exceptional.
In order to answer the question in general, we need pay particular attention to the topology we put on the loop space \cite{pacific}.

By abuse of terminology, a ``simplicial complex" refers to either the complex or its realization
(the underlying ``polyhedral space"). Let $X$ be a finite simplicial complex.
Endow $\Omega_{pl} (X)$, the set of all \textit{piecewise linear} loops on $X$, with a``weak topology" finer than the compact-open topology. This is a topology obtained as the colimit topology with respect to a canonical filtration of the piecewise linear PL loops (see \S\ref{plloops}). Let $\omega (X)$ be the corresponding space of based PL thin loops with the induced quotient weak topology. We prove

\bth\label{main} Let $X$ be a finite connected simplicial complex.
Then the space of thin PL loops $\omega (X)$ is a topological group homotopy equivalent to $\Omega (X)$, the space of based continuous loops.
\end{theorem}

This theorem is saying that at least for this particular family of loops, and with this topology, the thin fundamental group is in fact a group version of the loop space. Throughout this paper, we will write $\pi^1_{1}(X)$ for the thin loop space with the quotient topology induced from the compact open topology on $\Omega_{pl} (X)$, and we write $\omega (X)$ for the thin loop space which is a quotient of
$\Omega_{pl}^w(X)$, the space of PL loops with the weak topology.
The following diagram of continuous maps and homotopy equivalences $\simeq$ clarifies the relationship and gives a precise overview of the main results of the paper
$$\xymatrix{
\Omega^w_{pl} (X)\ar[r]^\simeq\ar[d]^\simeq& \Omega_{pl} (X)\ar[d]\ar[r]^\simeq&\Omega (X)\\
\omega (X)\ar[r]& \pi_{1}^1(X)
}$$
The vertical maps are quotient maps and the bottom map is a group homomorphism.
The necessity of departing from the compact-open topology on the righthand side, and considering the weak topology on the lefthand side, turns out to be a subtle point, one troubling issue being that the arclength functional (measuring the lengths of loops) is not continuous in the compact-open topology, in both the PL or PS categories, but it is in the weak topology. It can easily be checked that $\Omega_{pl}^w(X)$ is homotopy equivalent to $\Omega_{pl}(X)$ (see \S\ref{proofs}),
so by Theorem \ref{main} all spaces in the diagram except for $\pi_1^1(X)$, are homotopy equivalent. We are not able to check that $\omega (X)\simeq\pi_{1}^{1}(X)$ (open question). We must add that similar mapping space topology issues have risen in work of K. Teleman \cite{teleman}.

In \S\ref{loops} we discuss and summarize useful properties of various categories of loop spaces. In \S\ref{thin} and \S\ref{retracts} we introduce thin loop spaces and the weak topology. In \S\ref{plloops} we give a construction of $\omega (X)$ as a ``combinatorial" identification space and identify its homeomorphism type. In \S\ref{proofs} we prove our main result. Finally  in \S\ref{higher}, we define higher thin loop spaces $\pi_n^1(X)$, for $n\geq 1$, which consist similarly of all continuous loops in $\Omega^n(X)$, up to thin homotopy. In this case however one can check that as soon as $n>1$, $\pi_n^1(X)$ and $\Omega^n(X)$ are not generally of the same homotopy type.

Historically, Theorem \ref{main} has been approached by P. Gajer for the space of smooth loops. Based on a construction of Lefshetz, \cite{gajer} introduced a group model $G(X)$ for loop spaces when $X$ is a smooth manifold. This model is a quotient of $\Omega (X)$ by similar identifications as in the thin loop construction. In \cite{gajer} (Theorem 1.2) it is asserted that $G(X)$ is weakly homotopy equivalent to $\Omega (X)$ when $X$ is smooth. The proof is based on showing that $G(X)$ is the fiber of a universal principal bundle over $X$ with total space $E(X)$; the space of \textit{thin paths}. We believe that the proof of the contractibility of $E(X)$ is  incomplete and it is not clear to us how it can be filled in. One needs a more subtle approach to contractibility which is why our argument only works for PL loops and is in fact modeled over an earlier construction by Milnor \cite{milnor}.

Thin loop spaces are the main players in Teleman's work \cite{teleman}, although their homotopy type has not been addressed. They also appear ``in disguise" in Smale's thesis \cite{smale}. Smale considers regular curves on a Riemannian manifold $M$ which are classes of parameterized loops $\gamma$ such that $\gamma'(t)$ exists, is continuous and has non-zero magnitude for every $t\in I$. The equivalence relation is that $f\sim g$ if $f= g\circ h$ where $h$ is a diffeomorphism of $I$ with $h'(t)>0$ for every $t$.
It turns out that these loops are thin (\S\ref{thin}). We can write
$\omega_{reg}X$ for Smale's space of regular loops. Smale's main result (\cite{smale}, Theorem C) asserts that if $X$ is a Riemannian manifold,  there is a weak equivalence
$$\xymatrix{\omega_{reg}X \ar[r]^{w}&\Omega \tau X}$$
where $\tau X$ is the unit sphere bundle of $X$.
This paper suggests that if we allow piecewise regular curves, meaning that if we allow a finite number of points where the derivative is zero or undefined, up to suitable reparameterization, then the corresponding space of thin loops is weakly homotopy equivalent to $\Omega (X)$ but not to $\Omega\tau X$ (this we state as a second open question).

We conclude this paper by looking at the Borel construction associated to the conjugation action as in \cite{bc}, and deduce a cellular model for the thin free piecewise linear loops on $X$; $\ell (X)$, with the weak topology and show it is homotopy equivalent to the space of free loops $LX$ (Proposition \ref{freeloop}).

\vskip 5pt

{\sc Acknowledgments}: The ideas of this project were initiated at LICMAA'2015 in Beirut, Lebanon. We are much grateful to Roger Picken for sharing useful remarks on an early draft of the paper. We also thank Ali Maalaoui for discussions around the topic.


\section{Loops on Simplicial Complexes}\label{loops}

Throughout the paper, $X$ is a path-connected, finite simplicial complex thus compact.
A \textit{curve} means a parameterized path or loop.

Let $\Omega^n X$ be the space of basepoint preserving continuous loops from $S^n\subset\bbr^{n+1}$ to $X$.
This space is endowed with the compact-open topology which coincides for compact regular metric spaces
with the topology of uniform convergence with respect to the distance
function $sup_{x\in S^n}d(f(x),g(x))$ (\cite{hu}, Chapter VI, Lemma 2.2).
When $X$ is a smooth manifold, we denote by $\Omega^{n}_{sm}(X)$ the subspace of smooth loops.

For $X$ a compact smooth manifold, it is well-known that there is a homotopy equivalence $\Omega_{sm}^n(X)\simeq\Omega^n (X)$. Generally, for $N$ and $X$ smooth \textit{compact} manifolds, $\map{}^{sm}(N,X)$ (smooth maps) and $\map{}(N,X)$ (continuous maps) have the same homotopy type (see for example \cite{hansen} and references therein. See also Lemma \ref{pspl} below). In this paper we focus on piecewise linear (PL) and piecewise smooth (PS) loops on finite simplicial complexes. We will have to agree on what these are.

By definition, a finite  simplicial complex
is a space $X$ with a given finite triangulation.

\bde\label{defpl} Let $X$ to be a finite simplicial connected complex
with a chosen embedding $\phi$ in $\bbr^N$ so that all faces of $X$ are affine. This can always be arranged,
and we sometimes refer to $X$ as being a ``polyhedral space".
A path $\gamma : I\lrar X$ is PL, if the composite
$$\xymatrix{I\ar[r]^\gamma& X\ar@{^(->}[r]^{\phi} & \bbr^{N}}$$
is PL. A path is PS (piecewise smooth) if the same composite is PS.
\ede

\begin{properties}\label{remarks}.\rm
\begin{itemize}
\item 
The space $\Omega_{pl}X$ is topologized as a subspace of $\Omega (X)$.
As is indicated in Lemma \ref{pspl}, the homotopy type of $\Omega_{pl}X$ doesn't depend on the triangulation, nor does it depend on the embedding $X\hookrightarrow \bbr^N$.
\item A PL map $\gamma : I\rightarrow X$ has an underlying subdivision
$0=t_0< t_1< \cdots < t_n=1$ (not necessarily unique) such that the restriction
$\gamma_|: [t_i,t_{i+1}]\lrar\sigma\subset\bbr^{N}$ is linear, $\sigma$ being a face of $X$.
Given a subdivision as above for $\gamma$, it would be convenient to call $\gamma (t_i)=x_i$ a vertex of $\gamma$.
\item
If $\gamma$ is a PL loop (or path) with underlying subdivision $0=t_{0}< t_{1}< \dots < t_{n}=1$ as above, then its velocity is constant on each interval $]t_{i-1},t_{i}[$ and it is given by
$$ \| \gamma '(t)\|=\frac{d(x_{i},x_{i-1})}{t_{i}-t_{i-1}}, \quad t\in ]t_{i-1},t_{i}[$$
\item A loop $I\rightarrow X$ into a Riemannian manifold $X$ is \textit{regular} if it is smooth and has a non-zero derivative at every $t\in I$.
\item A smooth loop $\gamma : I\rightarrow X$
is \textit{uniform} if the velocity $ \| \gamma '(t)\|$ is the constant function. In particular
    $$\{\hbox{Non-constant Uniform}\} \hookrightarrow\{\hbox{Regular}\}$$
\item  A piecewise regular (PR) loop in a smooth manifold $X$ is a loop $I\rightarrow X$ with a subdivision $0=t_0<t_1 <\cdots <t_n=1$ such that $\gamma$ restricted to $[t_i,t_{i+1}]$, for any $i$, is either regular or constant (see \cite{smale},\S9, and \cite{teleman}, chapter III).  Piecewise regular loops can be defined on a finite simplicial complex embedded in $\bbr^N$ for some $N$, as for PL maps. PL loops or paths are also PR.
\item A  loop is PL uniform if $\gamma$ is either constant or $ \| \gamma'(t)\|=L>0$ at all but a finite number of points. This constant $L$ is forcibly the length of the curve.
\item In particular, a loop or path $\gamma$ is PL uniform iff it is constant or
$$t_{i}-t_{i-1}=\frac{d(x_{i},x_{i-1})}{\sum^{n}_{k=1}d(x_{k},x_{k-1})}  \quad 1\leq i\leq n $$
where again $x_i=\gamma (t_i)$.
\item When $\gamma$ is PL uniform, we say $t_i$ is redundant if $\gamma (t_i)\in [\gamma (t_{i-1}),\gamma (t_{i+1})]$. In this case $t_i$ can be suppressed from the subdivision without changing $\gamma$, and the subdivision gets reduced. Two subdivisions that vary by redundant vertices give the same PL uniform map. More on this in \S\ref{plloops}.
\item Composition of loops is not always well-defined in $\Omega_{sm}(X)$ since the composite is not always smooth at the basepoint. The composition is however well-defined for all piecewise (linear, regular, smooth, etc) functions.
\end{itemize}
\end{properties}

Let $P_{ps}X$ and $P_{pl}(X)$ be respectively the spaces of PS and PL paths
$I\lrar X$ sending $0$ to a fixed basepoint in $X$. We set throughout $\bullet$ to be PS or PL depending.
These spaces are contractible via the contraction pulling back a loop $I\times P_\bullet(X)\lrar P_\bullet(X)$ to the basepoint,
$(s,\gamma)\longmapsto \gamma_s$ where $\gamma_s(t) = \gamma ((1-s)t)$.
There is an evaluation map $ev: P_{\bullet}X\lrar X$, $\gamma\longmapsto \gamma (1)$, and the fiber over $x\in X$ is a copy of $\Omega_{\bullet} X$.

\bpr\label{evalfib} The evaluation maps $ev: P_{ps}X\rightarrow X$
and $ev: P_{pl}X\rightarrow X$
are quasifibrations.
\epr

\begin{proof} In a fixed triangulation, we can write $X$ as a union of simplexes $\bigcup_\alpha\sigma_\alpha$.
A lifting function for a map $p: E\rightarrow B$ assigns continuously
to each point $e\in E$ and path
$\gamma$ in $B$ starting at $p(e)$ a path $\lambda (e,\gamma )$ in $E$ starting at $e$ that is a lift for $\gamma$.
A map $p: E\rightarrow B$ is a fibration if and only if there exists a lifting function
for $p$ (Spanier, Theorem 8, Chapter 2).
Let's first of all show that for each simplex $\sigma$ of $X$, the restriction
$ev_|: ev^{-1}\sigma\lrar\sigma$ has a lifting function, and thus is a fibration. We will write it down for PL-paths, the proof for PS being totally similar. Pictorially this can be done as indicated in figure \ref{lifting1}. The simplex $\sigma$ is indicated by the solid simplex in the figure, and $\gamma$ is any path in $\sigma$ starting at $x$. We have to lift this path once we are given an element
of the preimage $e\in ev^{-1}(\sigma )$ which is a path ending at $x$, i.e. $ev(e)=x$.
The initial path is drawn as $e$ in the figure, ending at $x$. Given $e$, the lift of $\gamma$ is now given by the family of piecewise linear paths $e_t$ starting at $e(0)$ and ending at $\gamma(t)$ by extending linearly $e$ from
$x$ to $\gamma (t)$ as in the righthand side of Figure \ref{lifting1}. This linear extension can be done because
we are in $\sigma$.

\begin{figure}[htb]
\begin{center}
\epsfig{file=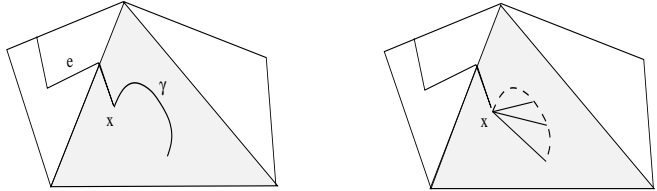,height=1in,width=3.4in,angle=0.0}
\caption{The lifting function for PL loops}\label{lifting1}
\end{center}
\end{figure}

The restriction of $ev$ over any closed simplex is therefore a fibration. We will show it
is a quasifibration over all of $X$ by (double) induction on cells of skeleta of $X$.
 The statement is trivial on the
zero skeleton $X^{(0)}$, so we suppose it is true on $L:=X^{(n)}$.
Let $\sigma$ be a simplex of dimension $n+1$. Let $a\in\sigma$ be an interior point. Then
$\sigma-\{a\}$ deformation retracts onto its boundary, and $L\cup(\sigma-\{a\})$ deformation retracts onto $L$ via a homotopy $d_t$. This deformation retraction is covered by a deformation retraction $D_t$ of
$\hbox{ev}^{-1}(L\cup(\sigma-\{a\}))$ onto $ev^{-1}(L)$, $D_0=id$, Image$(D_1)=\hbox{ev}^{-1} (L)$,
and $D_1: \hbox{ev}^{-1}(x)\lrar \hbox{ev}^{-1}(d_1(x))$ is a homotopy equivalence. By Lemma 1.3 of \cite{dold}, ev is now a quasibration over $L\cup (\sigma -\{a\})$. This means that $\hbox{ev}$ is  a quasifibration over each of the terms $L\cup (\sigma -\{a\})$ and $\mathring{\sigma}$ which are both open in $L\cup\sigma$.
By Lemma 1.4 of \cite{dold} again, $\hbox{ev}$ must be a quasifibration over $L\cup\sigma$. By induction it is a quasifibration over the $n+1$-skeleton and thus on all of $X$ (Lemma 1.5 in \cite{dold} says the conclusion is true even if $X$ is not necessarily finite).
\end{proof}

An immediate corollary to Proposition \ref{evalfib} is that for $X$ a finite simplicial complex, both $\Omega_{pl}(X)$ and $\Omega_{ps}(X)$ have the weak homotopy type of $\Omega (X)$. We can say more.
It is convenient to write hCW the category of spaces of the homotopy type of a CW complex.

\ble\label{pspl}
Both $\Omega_{pl}(X)$ and $\Omega_{ps}(X)$ are in hCW. In particular, both spaces are of the homotopy type of $\Omega (X)$.
\ele

\begin{proof}
Wlog, $X$ is a polyhedral space in $\bbr^n$ for some $n$.
It has an open tubular neighborhood $V$ that PL deformation retracts onto it.
We can then replace $\Omega_{\bullet}(X)$ by $\Omega_{\bullet}(V)$ up to homotopy, where $\bullet$ stands for either PL or PS. By a result of J.H.C.Whitehead, metrizable ANRs are in hCW (\cite{fp}, Theorem 5.2.3). Since both spaces are evidently metrizable (via the sup metric), one needs verify they are
ANRs.  Now $\Omega_\bullet (V)$ is an open subspace of $\Omega_\bullet(\bbr^n)$. An open subspace of an ANR is an ANR, so it is enough to show that $\Omega_\bullet(\bbr^n)$ is an ANR. But $\Omega_\bullet(\bbr^n)$ is a normed vector space and is locally convex. According to \cite{palais}, Theorem 5, it is an ANR. This shows indeed that $\Omega_\bullet (X)$ is in hCW. The quasifibrations in the earlier Proposition \ref{evalfib} are principal (contractible total spaces), so that $\Omega_\bullet X$ are weakly $\Omega (X)$. By being in hCW, they must be of the same homotopy type.
\end{proof}

Similarly, both Proposition \ref{evalfib} and Lemma \ref{pspl} have analogs for free loop spaces. A free loop is a map $[0,1]\rightarrow X$ with $\gamma (0)=\gamma (1)$. The free loop space is denoted by $LX$. We will write similarly $L_{ps}X$ and $L_{pl}X$ the spaces
of PS and PL free loops

\bpr\label{freeloopquasi} With $X$ a simplicial complex,
the evaluation maps $ev: L_\bullet X\rightarrow X$, $ev (\gamma)=\gamma (x_0)=\gamma (x_1)$ are quasifibrations for $\bullet=$PL and PL. As a consequence, $L_{ps}X\simeq LX\simeq L_{pl}X$.
\epr

\begin{proof} As in the proof of
Proposition \ref{evalfib}, the argument boils down to constructing a lifting function over a simplex $\sigma\subset X$; all the other details being the same.
To describe this lifting function,
let's look again at the righthand side of figure \ref{lifting1} where now
$e$ is a loop in $\pi^{-1}(\sigma)$ ending and starting at $x\in\sigma$. As in that
figure, $\gamma$ is a path in $\sigma$ starting at $x$. A lift of $\gamma(t)$ is the
loop starting at $\gamma (t)$, going linearly to $x$ in $\sigma$, traversing $e$ and then back to $\gamma (t)$ linearly. The remaining details are the same
and $ev$ is a quasifibration. The quasifibration $L_{ps}X\rightarrow X$ maps
to the loop space fibration $LX\rightarrow X$, and since both fiber and base are homotopy equivalent (Lemma \ref{pspl}), the total spaces must be weakly equivalent, thus equivalent since in hCW. Same for $L_{pl}X\simeq LX$.
\end{proof}


\section{Thin Loop Spaces}\label{thin}

We now define the notion of thin homotopy on spaces of loops. In this section we still assume the compact-open topology.

\bde\label{defthin} \cite{barrett}
A loop $\gamma$ is a \textit{thin loop} if there exists a homotopy (rel. the basepoint) of
$\gamma$ to the trivial loop with the image of the homotopy lying
entirely within the image of $\gamma$.  Two loops $f$ and
$g$ are said to be ``thin homotopic'' if there is a sequence of loops
$\gamma_1=f,\gamma_2,\ldots, \gamma_n=g$ such that for each $i\geq 1$,
$\gamma_i\circ \gamma_{i+1}^{-1}$ is thin.  Being thin homotopic defines an equivalence
relation on the space of based loops.
 \ede

\bde\label{repara} A ``reparameterization" of a path $\gamma : I\lrar X$ is any map obtained from $\gamma$ by precomposition with a continuous map $\alpha: I\lrar I$, $\alpha (0)=0, \alpha (1)=1$; i.e.
$$I\fract{\alpha}{\lrar} I\fract{\gamma}{\lrar} X$$
Depending on the category we work in (continuous, smooth, piecewise smooth or linear, etc), the reparameterization $\alpha$ will be a morphism in that same category.
The space of all continuous reparameterizations of $I$ fixing the endpoints is contractible. The space of piecewise smooth reparameterizations is also contractible, and in fact the space of diffeomorphisms of $I$ isotopically trivial and fixing the boundary; $Diff_0(I,\partial I)$, is trivial (see \cite{mann}, section 2).
\ede

\bex\label{killflare} Any two maps that differ by a reparameterization are thin homotopic.
Similarly, a ``flare" in a loop $\gamma$ is a path section that backtracks on itself. To be more precise, suppose this path section is uniformly parameterized.
We say $\gamma$ is a flare on $[-\epsilon,\epsilon]$ if it has the form:
$$\lambda : [-\epsilon, \epsilon]\lrar X\ \ ,\ \ \gamma (t) = \gamma (-t)\ ,\
0\leq t\leq\epsilon$$
This uniform flare is thin homotopic to the constant path
$c: [-\epsilon, \epsilon]\rightarrow X$, $c(t) = \gamma (-\epsilon)=\gamma (\epsilon)$ via the homotopy (in $s$)
$$H_s(t) = \begin{cases}\lambda \left((1-s)(t+\epsilon)-\epsilon\right) &, \ \ -\epsilon\leq t\leq 0\\
 \lambda ((1-s)(t-\epsilon) + \epsilon )&,\ \  0\leq t \leq \epsilon\end{cases}$$
 We say that $\gamma$ is a flare on $[a,b]$ if the composition
 $[-\epsilon,\epsilon]\rightarrow [a,b]\rightarrow X$ is a flare as already discussed, here
 the map $[-\epsilon,\epsilon]\rightarrow [a,b]$ is the linear map sending
 $(t-1)\epsilon + t\epsilon$ to $(1-t)a+tb$, $t\in [0,1]$.
\eex

\bex Figure 1 below shows two thin homotopic loops. In fact for $X=\Gamma$ a graph, any two homotopic loops in $\Gamma$ must be thin homotopic.
\eex

 \begin{figure}[htb]
\begin{center}
\epsfig{file=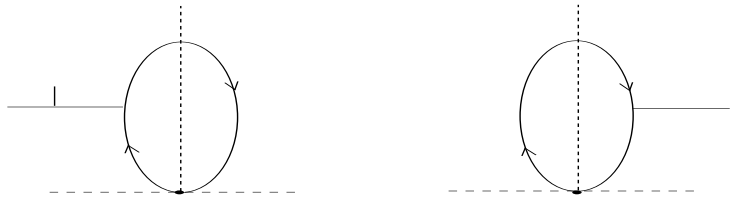,height=1in,width=3.8in,angle=0.0}
\caption{Two thin homotopic loops. Each loop has ``flares". Here a
\textit{flare} is a path section that goes back and forth once on itself.}
\end{center}
\end{figure}

For convenience we call a curve any parameterized path or loop.
Any two curves that differ by a reparameterization are thin homotopic.

Define the \textit{thin loop space} to be the quotient space
$${\pi_1^1}(X) = \Omega (X)/_\sim$$
where $\sim$ is thin homotopy. On $\pi_1^1 (X)$ there is defined a product
pairing by $[\gamma_1]*[\gamma_2]:=[\gamma_1\circ \gamma_2]$. This is well-defined.
Since reparameterization is obtained
by thin homotopy, $\pi_1^1 (X)$ is associative and has an identity which is the class of the constant loop at the basepoint. Moreover since for any loop
$\gamma$, the product $\gamma^{-1}\circ \gamma$ is a thin loop, its thin class $[\gamma^{-1}\circ \gamma]$ is trivial. But since $[\gamma^{-1}\circ \gamma]=[\gamma^{-1}]*[\gamma]$, it follows that $\pi_1^1 (X)$ has strict inverses $[\gamma ]^{-1}:=[\gamma^{-1}]$,
and $\pi_1^1 (X)$ is a group.  This group sits
between the loop space and the fundamental group as indicated in \eqref{sequence}. It is topologized as a quotient.

One can define thin loop spaces $\pi_1^1$ in any category of maps: R, PR, PS, PL (see Properties \ref{remarks}). A first needed condition is that the reparameterizations $\alpha : I\rightarrow I$ must be in the same category (see Definition \ref{repara}).  The second condition is that the homotopy between curves be adapted to the category.
In Definition \ref{pldef} we make precise the situation for PL loops which are our main focus.
Note that in \cite{barrett}, the thin loop space is for PS loops and
is a quotient by homotopies which are continuous and piecewise smooth on some ``paving" of $[0,1]^2$ consisting of polygons. On the other hand, in \cite{smale}, the regular curves considered on a Riemannian manifold $M$ are equivalence classes of standard regular curves under the relation that $f\sim g$ if $f= g\circ h$ where $h$ is a diffeomorphism of $I$ with $h'(t)>0$ for every $t$. In this case flares are automatically excluded and only reparameterizations with continuous and positive first derivatives are allowed. This means that a class of a Smale regular loop is the same as its thin regular homotopy class.

\bre (functoriality)
It is clear that a map
$f: X\rightarrow Y$ (within the same category) induces a homomorphism of groups $\pi_1^1 (X)\rightarrow\pi_1^1 (Y)$ (as in $\pi_1$), while an inclusion of spaces $X\subset Y$ induces a subgroup inclusion $\pi_1^1 (X)\subset\pi_1^1 (Y)$ (unlike $\pi_1$).
\ere

We will restrict in this paper to triangulated $X$ and PS or PL curves.
This allows us to use Lemma \ref{core} next.

\ble\label{core} Any piecewise regular curve on a simplicial complex $X$ can be
reparameterized into a piecewise uniform curve.
\ele

\begin{proof}
If the curve is constant there is nothing to prove. Otherwise,
the curve can be parameterized by arclength, meaning there is a parameterization
of the curve of the form
$[0,L]\rightarrow X$ with $\| {d\gamma\over dt}\| = 1$. This is proven in
 (\cite{tapp}, Proposition 1.25) but the same proof applies in the piecewise case
(see Lemma \ref{uniformretract}). Let $L$ be the length of the curve.
 By precomposing with the linear map $[0,1]\rightarrow [0,L], t\mapsto tL$, we obtain a uniform parameterization
$[0,1]\rightarrow X$ with $\| {d\gamma\over dt}\| ={L}$ at all but a finite
number of points.
\end{proof}

\bde\label{coreconst}
Consider a non-constant PR loop. It has only a finite number of flares by the very definition of PR (Properties \ref{remarks}).
 The \textit{core} of a piecewise regular curve $\gamma$ is the curve obtained as follows. First reparameterize $\gamma$ by arclength. If there are any flares, one uses the homotopy in Example \ref{killflare} to remove them one by one (order the flares by saying that a flare on $[a,b]$ comes before a flare on $[c,d]$ if $b\leq c$). With what is left, reparameterize again by arclength and remove any new flares. This process stops after a finite number of steps, and the end result is a parameterized loop with no flares. Take its uniform parametrization. This is called the ``core" of $\gamma$.
\ede

\bre Observe that the class of a thin loop depends not only on the shape of the image of that loop but also on the way
it is traversed. The core shall not be confused with the image of the loop.
\ere

The core construction above gives immediately the following characterization of thin homotopic
loops.

\bpr\label{criteria} Working with piecewise regular curves,
the following  are equivalent:\\\
(i) $\gamma_1$ and $\gamma_2$ are two thin homotopic loops\\
(ii) $\gamma_1$ and $\gamma_2$ have the same core (i.e they differ by a reparametrization or by addition and removal of flares).
\epr

\begin{figure}[htb]
\begin{center}
\epsfig{file=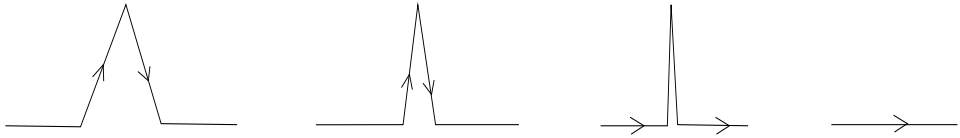,height=0.7in,width=5in,angle=0.0}
\caption{Topology in the space of thin loops: when both segments coincide, they ``vanish".}\label{touchvanish}
\end{center}
\end{figure}

Interestingly, the core construction is ``not continuous". In other words there is no section to the quotient projection $\Omega_\bullet X\lrar\pi_1^1 X$ (with $\bullet$ being PL or PR).
This is due to the following phenomenon summarized in figure \ref{touchvanish}.
One can try to see the lack of continuity in constructing a section $s: \pi_1^1 X\rightarrow \Omega_\bullet X$ by
arguing that because the composite
$$\pi_1^1 (X)\fract{c}{\lrar} \Omega_\bullet (X)\lrar \bbr^+\ \ ,\ \
[\gamma ]\longmapsto \hbox{core of $\gamma$}\longmapsto \hbox{length of core}$$
is not continuous, $c$ cannot be continuous. There is a problem in this argument however since the function that measures the length of piecewise smooth loops on a Riemannian manifold is not continuous in the $C^0$-topology; it is only upper-semicontinuous. It is however continuous in the $C^1$-topology as indicated in Lemma \ref{uniformretract}. We give below a more sound argument.

\ble\label{notcontinuous}
The ``core map" $c: \pi_1^1(X)\lrar\Omega_\bullet (X)$ sending a thin loop class to the core in that class is not a continuous map, thus it is not a section to the projection $q: \Omega_\bullet (X)\lrar\pi_1^1 (X)$.
 \ele

 \begin{proof} Set $h=c\circ q$, $X$ be Euclidean space and take a small open neighborhood $V$ of the
constant loop. Choose a loop $l$  which is not in $V$ but thin homotopic to the constant loop. Note that
$\ell\in h^{-1}(V)$.
Choose a sequence $(s_n)$ of loops which converges uniformly to $\ell$ and such that $h(s_n)$ is not in $V$. Then $(s_n)$ converges uniformly to a point in $h^{-1}(V)$ and $s_n$ is never in $V$. This means that
$h^{-1}(V)$ cannot be open and $h$, thus $c$, are not continuous.
\end{proof}

In this paper, we adopt the notation $\omega (X)$ for $\pi_1^1(X)$ (as a set) in the PL category, but with
a finer topology. More on this in \S\ref{thinmodel}.


\section{Retracts of Spaces of Loops}\label{retracts}

Throughout $X$ is a polyhedral space with a chosen basepoint $x_0$. We start with some properties of smooth loops.
Let $\Omega_{c^1}X$ be the space of regular curves (i.e. smooth with non-zero derivatives) in the $C^1$-topology. The subspace of uniform loops consists of
all loops $I\rightarrow X$ with $\| \gamma'(t)\| = {L}$, $L$ being the length of the curve $\gamma$ (which is always non-zero). The lemma below is stated for loops in the $C^1$-topology. Briefly, this is the initial topology defined by the map $T: \Omega_{c^1}(X)\rightarrow \Omega (X)\times Map(TS^1,TX), f\mapsto (f,df)$ (\cite{hirsch}, chapter 2.4). The lemma below is of independent interest but is not really used elsewhere.

\ble\label{uniformretract} Let $U\Omega_{c^1} (X)$ be the subspace of all uniform loops in $\Omega_{c^1}X$ in the $C^1$-topology. Then $U\Omega_{c^1}X$ is a deformation retract of $\Omega_{c^1}X$.
\ele

\begin{proof}  Let $\gamma : I=[0,1]\rightarrow X$ be a regular curve and $L_{\gamma}=\int_0^1||\gamma'(u)||du > 0$ its length.
Consider the normalized arc-length function
$$\theta_{\gamma}: I\rightarrow I,
\quad \theta_{\gamma}(t) = \frac{1}{L_{\gamma}}\int_0^t||\gamma'(u)||du$$
$\theta_{\gamma}$ is a well defined $C^1$-diffeomorphism. Therefore there is a unique well defined smooth curve
$\tilde\gamma: I\rightarrow X$ satisfying $\tilde\gamma(\theta_{\gamma}(t))=\gamma(t).$
We see that
$\tilde\gamma'(\theta_{\gamma}(t)) = \frac{L_{\gamma}\gamma'(t)}{||\gamma'(t)||}$. In particular  $||\tilde\gamma'(s)||=L_{\gamma}$  and $\tilde\gamma$ is uniform.

Let $i: U\Omega_{c^{1}}(X)\longrightarrow \Omega_{c^{1}}(X)$
be the inclusion map and $r: \Omega_{c^{1}}(X)\longrightarrow U\Omega_{c^{1}}(X)$
the map that takes a loop $\gamma$ to its uniformization $\tilde{\gamma}$ (we checked it is well-defined and this map is continuous in the $C^1$-topology).
Moreover $r\circ i=id_{U\Omega_{c^{1}}(X)}$, thus $r$ is a retraction. Define
$H: \Omega_{c^{1}}(X)\times I\longrightarrow \Omega_{c^{1}}(X)$ by
$$H(\gamma,s)(t)=\tilde\gamma((1-s)\theta_{\gamma}(t)+st)$$
Then $H$ is a strong deformation retract from $id_{\Omega_{c^{1}}(X)}$ to $i\circ r$.
\end{proof}

The proof above works since the length of curves is a continuous function in the $C^1$-topology. This is not the case in the $C^0$-topology. We can remedy the $C^0$ pathology by introducing a slightly different topology on the space of piecewise regular (or linear) loops.

\bde Fix a  simplicial complex $X$ (with given finite triangulation). Filter $\Omega_{pl}X$ by the subspaces
$F_n\Omega_{pl}X$ consisting for every $n$ of all piecewise loops which have a subdivision of $I$ with at most $n$ $t_i$'s, i.e. $t_0=0\leq t_1\leq\cdots\leq t_k= 1$ for $k\leq n$. These are loops which can be represented by at most $n$ segments.  The topology on $F_n\Omega_{pl}X$ is the sup-metric topology. We have a filtration
$\cdots\subset F_n\Omega_{pl}X\subset F_{n+1}\Omega_{pl}X\cdots$
and as sets
$$\Omega_{pl}X = \bigcup F_n\Omega_{pl}X$$
The weak topology is defined to be the colimit or final topology on $\Omega_{pl}X$
(\cite{spanier}, page4).
We write $\Omega_{pl}^wX$ for this space.
\ede

Let $U\Omega_{pl}X$ to be the subspace of piecewise uniform loops as in Properties \ref{remarks}, $i: U\Omega_{pl}X\hookrightarrow \Omega_{pl}X$ the inclusion, and $F_nU\Omega_{pl}X = F_n\Omega_{pl}X\cap U\Omega_{pl}X$ the corresponding filtration term.

\ble\label{firstlem}
$F_nU\Omega_{pl}X$ is a strong deformation retract of $F_n\Omega_{pl}X$. These retracts for all $n$ induce a deformation retract $r$ of the piecewise linear loops $\Omega^w_{pl}X$ onto the piecewise uniform loops $U\Omega^w_{pl}X$ (in the weak topology).
\ele

\begin{proof}
Given a PL loop or path $\gamma$ and $0=t_{0}\leq t_{1}\leq\dots\leq t_{n}= 1$ such that $\gamma$ is linear on each segment $[t_{i-1},t_{i}]$ and $x_{i}=\gamma(t_{i})$, $1\leq i\leq n$,
define a new subdivision $(\tilde{t}_{i})_{0\leq i \leq n}$ of the unit interval by $\tilde{t}_{0}=0$ and for $1\leq i\leq n$
\begin{eqnarray*}
\displaystyle\tilde{t}_{i}=
\begin{cases}\tilde{t}_{i-1}+\frac{d(x_{i},x_{i-1})}{\sum^{n}_{k=1}d(x_{k},x_{k-1})} ,& \hbox{if $\gamma$ is non-constant}\\
t_i,&\hbox{if $\gamma$ is constant}
\end{cases}
\end{eqnarray*}
Let $\tilde{\gamma}$ be the PL path which is linear on each segment $[\tilde{t}_{i-1},\tilde{t}_{i}]$ and satisfying  $\tilde{\gamma}(\tilde{t}_{i})=x_{i}$ (Properties \ref{remarks}). Then $\tilde{\gamma}$ is the uniform reparametrization of $\gamma$.
Note that there is no unique subdivision associated to the map $\gamma$. Indeed if in the subdivision $(t_i)_{0\leq i\leq n}$,
$\gamma$ is uniform on $[t_{i-1},t_{i+1}]$ and $x_i=\gamma (t_i)\in [x_{i-1},x_{i+1}]$ in some face of $X$, then
the new subdivision we get by omitting $t_i$ still represents $\gamma$. This point is discussed in more details in the next section
and in Lemma \ref{homeo}. In all cases, it can be shown that $\gamma$ does not depend on this representation, and that
the $\tilde{t}_i$ depend continuously on $\gamma$, and not again on the choice of subdivision.
The map
$$r_{n}: F_n \Omega_{pl}X\longrightarrow F_n U\Omega_{pl}X$$
which takes a loop $\gamma$ to its uniformization $\tilde{\gamma}$ is well-defined and continuous.
When $\gamma$ is uniform, $\tilde\gamma=\gamma$ and $r_n$ is a retract of the inclusion
$i_{n}: F_n U\Omega_{pl}X\longrightarrow F_n\Omega_{pl}X$.

For $\gamma$ and $\tilde{\gamma}=r_n(\gamma)$, $(t_{i})_{0\leq i \leq n}$, $(\tilde{t}_{i})_{0\leq i \leq n}$ and $(x_{i})_{0\leq i \leq n}$ as above, and $s\in [0,1]$, define $t^{s}_{i}=(1-s)t_i + s\tilde{t}_{i}, \quad 0\leq i \leq n.$ Define $H_n(\gamma,s)$ to be the PL loop in $X$ which is linear on $[t^{s}_{i-1},t^{s}_{i}]$ and satisfying $H_n(\gamma,s)(t_{i}^{s})=x_i.$ Then
$H_n$ is the desired strong deformation retract from $id$ to $i_n\circ r_n$.
The maps $r_n$ induce in turn a retraction
$r:  \Omega^w_{pl}X\longrightarrow  U\Omega^w_{pl}X$
and the homotopies $H_n$ induce a strong deformation retract $H$ from $id$ to $i\circ r$.
\end{proof}


\section{Piecewise Linear Loops: The Model}\label{plloops}

From now on, we restrict to PL loops $\gamma : I\rightarrow X$, $\gamma (0)=\gamma (1)$, with PL homotopies $F: I\times S^1\rightarrow X$ (meaning the restriction $F_t$ is PL for any $t\in I$) and PL-reparameterizations. This section will characterize (combinatorially) the thin PL loop spaces $\omega (X)$ following earlier ideas of Milnor \cite{milnor} (see appendix).

\begin{figure}[htb]
\begin{center}
\epsfig{file=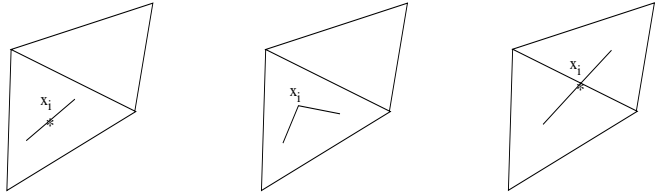,height=1in,width=3.4in,angle=0.0}
\caption{(a) $x_i$ is redundant. In (b) and (c) $x_i$ is essential}\label{redundancy}
\end{center}
\end{figure}

Let $\gamma : I\to X$  be a uniform PL loop as in Definition \ref{defpl}. This defines a tuple $(\gamma (t_0),\gamma (t_1),\ldots, \gamma (t_{n}))$ in $X^{n+1}$ with $\gamma (t_0) = x_0 = \gamma (t_{n})$ and $\gamma (t_i),\gamma (t_{i+1})$ in the same simplex. We say $x_i=\gamma (t_i)\in X\subset\bbr^n$ is a ``redundant point" if $x_{i}\in [x_{i-1},x_{i+1}]$ within the same simplex. Figure \ref{redundancy} illustrates redundancy.

Define
\begin{equation}\label{sk}
\tilde S_k = \{(x_0,x_1,\ldots, x_{k-1},x_k)\in X^{k+1}\ \ ,\ \
x_k=x_0\ ,\ x_i,x_{i+1}\ \hbox{in the same simplex of $X$}\}
\end{equation}
and define the quotient
$\displaystyle W_{n} := \coprod_{k=0}^{n} \tilde S_k/_\sim$
where $\sim$ is the identification relation determined by
$$(x_0,x_1,\ldots, x_k)\sim (x_0,x_1,\ldots, \hat x_i,\ldots, x_k)\ \ ,\ \ \hbox{if $x_i$ is redundant}$$
One then defines $W_\infty$ to be the direct limit $\bigcup W_n$ with the final topology.
The following lemma is equivalent to saying that a piecewise linear loop is uniquely defined by its non-redundant vertices.

\ble\label{homeo} $W_\infty$ is homeomorphic to $U\Omega^w_{pl}X$.
\ele

\begin{proof}
There is a map $\phi: W_{\infty}\rightarrow U\Omega_{pl}X$ which
to a class ${\bf x}=[x_0,x_1,\ldots, x_{n-1}, x_{0}]\in W_n$ associates the piecewise uniform linear loop
$\phi ({\bf x})=\gamma$ with subdivision $0=t_0\leq t_1 \leq\cdots \leq t_n= 1$ such that $\gamma (t_i)=x_i$. This map is well-defined since
a map that is uniformly linear on $[t_{i-1},t_i]$ and $[t_i,t_{i+1}]$ is identified to a uniform linear map on
$[t_{i-1},t_{i+1}]$ provided $x_i\in [x_{i-1},x_{i+1}]$.
The map $\phi$ is continuous on $W_\infty$ because its restriction to every $W_n$ is continuous.
Each
$\phi_n: W_{n}\rightarrow F_{n}U\Omega_{pl}X$ is clearly bijective, with inverse
$\psi_n: F_nU\Omega_{pl}X\lrar W_n$
sending a PL map $\gamma$ with subdivision $0=t_0\leq t_1\leq \cdots \leq t_n= 1$ to $\psi (\gamma ) = [x_0,\gamma (t_1),\ldots, \gamma (t_{n-1}), x_0 ]$.
Since
$W_{n}$ is compact and $F_{n}U\Omega_{pl}X$ is Hausdorff,  $\phi_n$ is a homeomorphism. It follows that the colimit map $W_{\infty}\rightarrow U\Omega^w_{pl}X$ is a homeomorphism.
Observe that all maps below
$$W_\infty\fract{}{\lrar}\Omega^w_{pl}X\fract{r}{\lrar} U\Omega^w_{pl}X\fract{\psi}{\lrar} W_\infty$$
are continuous, with $r$ as in Lemma \ref{firstlem}, and that the composite is the identity.
\end{proof}

\subsection{Thin PL-Loops}\label{thinmodel}
Two PL maps $I\lrar X$ are PL-homotopic if there is a PL-homotopy $I\times I\lrar X$, rel the basepoint,
between both maps (see Definition \ref{pldef} ahead).
We then define $\omega (X)$ to be the quotient of $\Omega^w_{pl}(X)$ as in
Definition \ref{defthin}. This space turns out to have a ``combinatorial description"
which is, up to homeomorphism, a quotient of $W_\infty$. We make this precise below.

\begin{figure}[htb]
\begin{center}
\epsfig{file=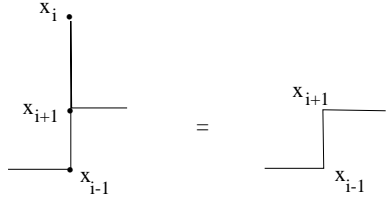,height=1.2in,width=2in,angle=0.0}
\caption{An identification $\approx$ of two ``combinatorial" thin loops
}\label{identification}
\end{center}
\end{figure}

\bde\label{milnorid} With $\tilde S_k$ as in \eqref{sk}, let $S_n(X)$ (or $S_n$ for short)
be the quotient space
\begin{equation*}
S_n(X) := \coprod_{0\leq k\leq n}\tilde S_k/_\approx
\end{equation*}
where $\approx$ is the equivalence relation identifying
\begin{equation}\label{identify}
(x_0,\ldots, x_{i-1}, x_i, \ldots, x_k)\approx (x_0,\ldots, x_{i-1},\hat{x}_i,x_{i+1},\ldots, x_k)\ \ \ ,\ \ \ x_k=x_0
\end{equation}
whenever  the triple $x_{i-1},x_i,x_{i+1}$ are aligned within the same simplex.
Here $\hat x_i$ means the $i$-th entry is deleted.
Figure \ref{identification} illustrates this identification which
is designed to kill flares (or backtracking).
We will write an equivalence class of an element in $S_n(X)$ as $[x_0,x_1,\ldots, x_n,x_0]$ and generally represent this tuple in reduced form (i.e. when no further identifications can be made).
Note that redundancy corresponds to when $x_i\in [x_{i-1}, x_{i+1}]$ so that
$S_n$ is a further quotient of $W_n$.
We define $\displaystyle S(X) = \bigcup_{n\geq 1}S_n(X) = \coprod_{k\geq 1}\tilde S_k/_\approx$.

\ede

\bth\label{pltos}
There is a homeomorphism $S(X)\cong\omega (X)$.
\end{theorem}

\begin{proof} Both maps $\psi$ and $\phi$ of Lemma \ref{homeo}
restrict after passing to quotients to maps $f$ and $g$ respectively as follows
$$\xymatrix{W_\infty\ar[r]^\phi\ar[d]&\Omega^w_{pl}(X)\ar[r]^\psi\ar[d]&W_\infty\ar[d]\\
S(X)\ar[r]^f&\omega(X)\ar[r]^g&S(X)}
$$
Vertical maps are (continuous) quotient maps, and both $\phi,\psi$ are continuous in the weak topology. By diagram chasing, both $g$ and $f$ are continuous. But $g$ is bijective, so $g=f^{-1}$ and both spaces are homeomorphic.
\end{proof}

Recall that $\omega (X)$ sits between $\Omega_{pl} X$ (an $H$-space) and
$\pi_1(X)$ (a group in the discrete topology).

\bco\label{top} $\omega (X)$ is a topological group.
\eco

\begin{proof} Both left and right translations are continuous and $\omega (X)$ is a semigroup. The space $\Omega^w X$ is compactly generated (or CG) being a colimit of a small diagram in CG (\cite{neil}, Prop. 2.23). It then follows that the product of quotient maps
$\Omega^w_{pl}X\times\Omega^w_{pl} X\lrar \omega (X)\times\omega (X)$ is also a quotient map (see \cite{neil}, Prop. 2.20). The commutativity of the diagram below implies that the composition $\ast$ below is continuous as long as loop sum above $\ast$ is continuous
$$\xymatrix{\Omega^w_{pl}X\times\Omega^w_{pl}X\ar[r]^{\ \ \ \ast} \ar[d]&\Omega^w_{pl}X\ar[d]\\
\omega (X)\times\omega (X)\ar[r]^{\ \ \ \ast}&\omega (X)}$$
The inverse map is also continuous and the proof is complete.
\end{proof}

The following two remarks describe what can happen with thin loops when working in the continuous and piecewise smooth categories.

\bre In general one can topologize the fundamental group as a quotient of $\Omega (X)$
(denoted by $\pi_1^{qtop}(X)$). Then each of the maps
$\Omega (X)\rightarrow\pi_1^{qtop}(X)$ is continuous, and
$\pi_1^{qtop}(X)$ is a semi-topological group. It is not true in general working in this category that $\pi_1^{qtop}(X)$ is a topological group as pointed out in \cite{bf}.
\ere

\bre\label{thepscase} (the piecewise smooth case).
Let $X$ be a smooth manifold and $\Omega_{ps}(X)$ its space of piecewise smooth loops.
Similarly as above, thin homotopy defines an equivalence relation on these loops and we have a thin loop group for this family of smooth loops. The thin loop group $\omega_{ps}X$ acts on the space of thin piecewise smooth paths $P^{thin}_{ps}(X)$ and the projection $P^{thin}_{ps}(X)\lrar X$, $\gamma\mapsto \gamma (1)$ has the structure of a principal bundle with group $\omega_{ps}X$.
The proof of this fact proceeds as in Lemma \ref{principal} with the distinction that we make use of the existence of a unique geodesic between two close enough points.
As pointed out in the introduction, it is not clear how to show this bundle is universal;
i.e. that $P^{thin}_{ps}(X)$ is contractible in this category.
\ere


\section{Proof of the Main Results}\label{proofs}

The idea of the proof of Theorem \ref{main} is standard and relies on constructing a contractible space of thin paths. The details for the contraction are surprisingly technical and are dealt with in this section.

Extend thin homotopy to the space of paths: we say two paths $\gamma_1,\gamma_2 : I\rightarrow X$, $\gamma_1(0)=\gamma_2(0)$, $\gamma_1(1)=\gamma_2(1)$,
are thin homotopic if the loop $\gamma_2\circ\gamma_1^{-1}$ is thinly homotopic to the constant loop at $\gamma_1(0)$. Two paths are thin homotopic if they have the same endpoints and differ by flares and reparameterizations. We will write $P^{thin}X$ the equivalence classes of
thin paths starting at $x_0\in X$.
Define
\begin{equation}\label{tildee}
\tilde E_k = \{(x_0,x_1,\ldots, x_k)\in X^{k+1} \ ,\ \ x_i,x_{i+1}\ \hbox{in the same simplex of $X$}\}
\end{equation}
This contains the subspace $\tilde S_k$ consisting of $x_{k}=x_0$ \eqref{sk}. Let $E_n$ be the quotient space
$$
E_n := \bigcup_{0\leq k\leq n}\tilde E_k/_\approx
$$
where $\approx$ is the same (thin) equivalence relation \eqref{identify}. We will at times write
$\approx$ as $\sim_{thin}$ for clarity. The direct limit construction
$E(X) = \bigcup E_n$, endowed with the colimit topology
gives a model for the thin path space (as in the case of thin loops in Theorem \ref{pltos}).
We will argue in this section that the $E_n$'s are contractible and that the colimit $E$ is contractible.

The starting point is a proposition that summarizes some of main results about simplicial complexes that we need. These are collected from the book \cite{rourke}. We give however a slightly different definition of PL maps in our case. Let $K\subset\bbr^n$ be a polyhedron.

\bde\label{pldef} A map $f:K\lrar \bbr^m$ is \textit{piecewise linear} if there exists a triangulation $\{\sigma_i\}_{i\in I}$ such that each of the restrictions $f|_{\sigma_i}$ is linear (that is the restriction of an affine map from $\bbr^n$ to $\bbr^m$).
If $K\subset\bbr^n$ and $L\subset\bbr^m$ are polyhedra, we say that a map $f: K\lrar L$ is piecewise
linear (or PL) if the underlying map $f: K\lrar\bbr^m$ is piecewise linear.
\ede

\bth\label{subdivisions}
Let $K$ and $L$ be abstract simplicial complexes (so that $|K|$ and $|L|$ are the underling
polyhedra).  We write $L'\triangleleft L$ if $L'$ is a sudivision of $L$. Then
\begin{enumerate}[(i)]
\item (\cite{rourke}, 2.16) Suppose that $f: K\rightarrow L$ is simplicial and $L'\triangleleft L$. Then there is $K'\triangleleft K$ such that $f: K'\rightarrow L'$ is simplicial. The dual of this statement is not true.
 \item (\cite{rourke}, 2.17) Suppose that $f: L\rightarrow K$ is a simplicial injection and $L'\triangleleft L$. Then there is $K'\triangleleft K$ such that $f: L'\rightarrow K'$ is simplicial.
\item (\cite{rourke}, 2.12) If $|L_i|\subset |K|$, $i=1,\ldots, r$, then there are simplicial subdivisions $K'\triangleleft K$ and $L'_i\triangleleft L_i$ such that $L'_i\subset K'$ for each $i$.
\item (\cite{rourke}, 2.14) Let $f: |K|\rightarrow |L|$ be PL, then there are sudivisions $K'\triangleleft K$ and $L'\triangleleft L$ such that $f: |K'|\rightarrow |L'|$ is simplicial\footnote{In fact in \cite{panov} this is used as a definition of PL; i.e. a PL map is a simplicial map between some subdivisions.}.
\end{enumerate}
\end{theorem}

\ble\label{cwstruct}
The inclusion $E_{n-1}\hookrightarrow E_n$ is a cofibration.
\ele

\begin{proof}
Define
\begin{equation}\label{degenerate}
D_n =\{(x_0,\ldots, x_n)\in S_n\ |\ x_{i-1},x_i,x_{i+1}\ \hbox{aligned in some simplex $\sigma\in X_*$ for some $i$}
\}
\end{equation}
the ``degenerate" subspace.
We have a quotient map $\beta_n: D_n\lrar E_{n-1}$ and a pushout diagram
 \begin{equation}\label{pushout}
 \xymatrix{
 D_n\ar[d]^{\beta_n}\ar@{^(->}[r]&\tilde E_n\ar[d]\\
 E_{n-1}\ar[r]&E_n}
 \end{equation}
where $\tilde E_n$ is as in \eqref{tildee}.
Let $X_*=\{\sigma_i\}_{i\in I}$ be a finite simplicial decomposition of $X$ with $x_0$ a chosen vertex. After putting an order on the vertices of this decomposition, we can construct in a standard way a simplicial decomposition for $X^n$, denoted $X^n_*$.
Now by construction, $\tilde E_n$ can be written as
$$\tilde E_n = \bigcap_{k\geq 2}\left(\bigcup_{i\in I} X^k\times\sigma_{i}^2\times X^{n-k-1}\right)
\cap \bigcup_{\tau_j\ni x_0}\left(x_0\times\tau_j\times X^{n-1}\right)$$
and this inherits the structure of a subcomplex of $X^n_*$ as an intersection of subcomplexes.

We first argue that $(\tilde E_n,D_n)$ is a simplicial pair (for some simplicial structure).
Note that if $\sigma$ is a cell, then there is a map
$\tau_{13}: \sigma\times I\times \sigma \lrar \sigma^3$
which sends $(x,t,y)\longmapsto (x, tx+(1-t)y, y)$. This is a homeomorphism onto its image as long as $x\neq y$. It is a PL map.
Similarly there is a map
$\tau_{12}: \sigma\times\sigma\times I\lrar \sigma^3$,
$(x,y,t)\longmapsto (x,y, tx+(1-t)y)$ and a third map $\tau_{23}$ so that
\begin{equation}\label{themap}
\tau_{12}\sqcup \tau_{13}\sqcup \tau_{23}: I\times\sigma\times\sigma\sqcup
\sigma\times I\times\sigma\sqcup \sigma\times\sigma\times I\lrar\sigma^3
\end{equation}
is a PL map. By Theorem \ref{subdivisions}-(iv), there are simplicial decompositions of
both spaces in \eqref{themap} so that the map becomes simplicial. As an immediate
consequence we get a simplicial structure on
$$\overline{\sigma} := Im(\tau_{12})\cup Im(\tau_{13})\cup Im(\tau_{23})\ \subset\ \sigma^3$$
The degenerate set can be written in turn as
$$D_n = \left(\bigcup_{i\in I} X^k\times\overline{\sigma}_{i}\times X^{n-k-2}\right)\cap \tilde E_n$$
Each term of the union is simplicial. This means that for appropriate subdivisions of $\tilde E_n$, the pair $(\tilde E_n,D_n)$ is simplicial (this is Theorem \ref{subdivisions} -(iii)). Since $D_n$ is subsimplicial in $\tilde E_n$, the inclusion is a cofibration. Since cofibrations are preserved under pushouts, the claim of the proposition becomes a consequence of \eqref{pushout}.
\end{proof}

\bre We believe that $(E_n,E_{n-1})$ can be given the structure of a CW pair,
and that the space of thin maps is a CW complex.
\ere

The following is a reminder of some standard  properties. All spaces below are assumed to be
compactly generated, of the homotopy type of CW complexes.

\ble\label{XoverA}
Suppose $A\hookrightarrow X$ is a cofibration.\\
(i) If both $A$ and $X/A$ are contractible, then $X$ is contractible.\\
(ii) If both $A$ and $X$ are contractible, then $A$ is a deformation retract of $X$.
\ele

\begin{proof}
The following pushout diagram
$$\xymatrix{A\ar[r]\ar[d]&X\ar[d]\\
\star\ar[r]&X/A}$$
is a homotopy pushout, with the top map being a cofibration. It follows that
when $A$ is contractible, and the left vertical map $A\rightarrow *$ is a homotopy equivalence, the right vertical map $X\simeq X/A$ is an equivalence as well.
Consequently $X/A$ is contractible if and only if $X$ is contractible.
On the other hand, $X/A$ contractible if and only if $A$ is a deformation retract of $X$.
\end{proof}

\bre If $A\hookrightarrow X$ is not a cofibration (in particular $(X,A)$ is not a CW pair), then it can happen that both $A$ and $X/A$ are contractible, but $X$ isn't. A famous counterexample is the ``Griffiths Twin Cone" which is a non-reduced cone on a wedge of two Hawaiian rings $X=H\vee H$ (here can choose $A$ to be one of the wedge copies of $H$) \cite{griffiths}.
\ere

\bpr\label{contractible} The space $E(X)$ is contractible.
\epr

\begin{proof} We start by showing that each $E_n$ is contractible.
We proceed inductively. Let $D_n\subset E_n$ be the ``degenerate" subspace defined in \eqref{degenerate}.
From the pushout \eqref{pushout}, we have a homeomorphism
$${E_n\over E_{n-1}}\cong {\tilde E_n\over D_n}$$
In $\tilde E_n$ one contracts $(x_0,\ldots, x_{n},x_{n+1})$ to
$(x_0,\ldots, x_{n-1},x_{n},x_n)$ via the standard contraction $F$ pulling $x_{n+1}$ to $x_{n}$ linearly.
The final effect of this homotopy is a map $F_1: \tilde E_n\rightarrow D_{n}$ while at every stage
$F_t: D_n\rightarrow D_n$. This homotopy induces a contraction of $\tilde E_n/D_n$ so
this quotient is contractible\footnote{Note that in the case of Milnor's construction, see \S\ref{milnormodel}, if we set the degenerate set there to be
$$
D_n^M =\{(x_0,\ldots, x_{n})\in \tilde E_n\ |\ x_{i-1}=x_{i+1}\ \hbox{for some $i$}
\}
$$
then the homotopy $F$ doesn't leave $D^M_n$ invariant. A different approach for contractibility is required.}. Observe that $E_n$ is of the homotopy type of a CW complex inductively by \eqref{pushout}, with the two top spaces being
already simplicial. We can now deduce the contractibility of $E_n$ by induction as well. Evidently $E_0=\{x_0\}$ is contractible.
Assume $E_{n-1}$ is contractible. Since
${E_n/ E_{n-1}}=\tilde E_n/D_n$ is contractible, and since $E_{n-1}$ is contractible, it follows by both Lemma \ref{cwstruct} and Lemma \ref{XoverA} that
$E_n$ is contractible, so by induction this is valid for all $n\geq 0$.

To show that $E(X)$ is contractible in the weak topology, we observe first that since
again $E_{n-1}$ into $E_n$ is a cofibration, and that both
$E_n$ and $E_{n-1}$ are contractible, then by Lemma \ref{XoverA},
$E_{n-1}$ is a deformation retract of $E_n$. The contraction of $E_n$ to a point can be chosen to be first the deformation retraction of $E_n$ onto $E_{n-1}$ followed by the
contraction of $E_{n-1}$ to point. This means that we can choose contractions compatibly so that the following commutes for each $n$
$$\xymatrix{
E_n\times I\ar[r]&E_n\ar[r]&E(X)\\
E_{n-1}\times I\ar@{^(->}[u] \ar[r]&E_{n-1}\ar@{^(->}[u]\ar[ur]&
}$$
By definition of the colimit and the weak topology, these compatible maps induce a continuous map
$E(X)\times I\lrar E(X)$ which is a contraction (i.e. a homotopy between the identity and the constant map).
\end{proof}

\bre The contractibility of $E_n$ was not so trivial to obtain for the following reason.
Consider
$$\bigcup\tilde E_k/\sim \ \ \ ,\ \ \
(x_0,\ldots, x_k)\simeq (x_0,\ldots, \hat{x}_i,\ldots, x_k)\ \hbox{if}\ x_i=x_{i+1}$$
In \eqref{tildee}, there is a ``contraction" to basepoint which
takes a tuple $(x_0,\ldots, x_k)$ and moves consecutively each last entry $x_k$ to $x_{k-1}$ along the segment $[x_{k-1},x_{k}]$.
This contraction doesn't descend to $E_n$ for the following reasons. If it is written as such (\cite{gajer}, page 206), then this presupposes the existence of a continuous choice of a representative or a section $E_n\rightarrow \tilde E_n$ which we don't believe
exists according to the ``core Lemma" \ref{notcontinuous}.
If one writes the retraction on representatives, thus as a map
$I\times E_n\lrar E_n$, then it won't be well-defined.
\ere

We now turn to $S(X)\subset E(X)$, the subspace consisting of elements (or loops)$[x_0,x_1,\ldots, x_n,x_0]$, for some $n$.
There is a well-defined action of $S(X)$ on $E$ given by
$$[x_0,x_1,\ldots, x_{n},x_0]\cdot [x_0,y_1,\ldots, y_{m}] \longrightarrow [x_0,\ldots, x_{n}, x_0, y_1,\ldots, y_{m}]$$
which restricts to a group structure on $S(X)$ having inverses
$$[x_0,x_1,\ldots, x_{n},x_0]^{-1}= [x_0,x_{n},\ldots, x_1,x_0]$$
The group $S(X)$ acts freely on $E(X)$.

\ble\label{principal}
The triple $(S(X),E(X), X)$, with projection $\pi: [x_0,\ldots, x_n]\longmapsto x_n$,
is a principal fiber bundle.
\ele

\begin{proof} The proof runs similar to that of Milnor.
We recall a \textit{principal} $G$-bundle is a topological space $P$ with a continuous free action of $G$ such that the projection $\pi: P\lrar P/G$ is locally trivialized (for space and action); i.e. for every $x\in P/G$, there is an open neighborhood $U$ so that $\pi^{-1}(U)\cong G\times U$ and the action of $G$ on $\pi^{-1}(U)$ corresponds to $g(h,x)=(gh,x)$. Here we claim that $P=\tilde{E}$ has a principal action by $G=S(X)$.
The key point is the local triviality. For every $x\in X$,
pick $U$ to be the star neighborhood of $x$ in $X$, and choose
$[x_0,x_1,\ldots, x_{n-1},x]$ a fixed path ending in $x$, i.e. an element of $\pi^{-1}(x)$.
Then the map
$$\phi_x : G\times U\lrar \pi^{-1}(U)\ \ ,\ \ (g,y)\longmapsto g\cdot [x_0,x_1,\ldots, x_{n-1},x]\cdot [x,y]$$
is a homeomorphism with inverse
$$[x_0,y_1,\ldots, y_{n-1}, y]\longmapsto ([x_0,y_1,\ldots, y_{n-1},y, x,x_{n-1},\ldots, x_1,x_0],y)$$
One point here is the continuity of both the map and its inverse. The second point is that the thin relations $\approx$ ensure this is a bijection.

The action $S(X)\times E(X)\rightarrow E(X)$ is given by concatenation at the basepoint (or composition)
\begin{equation}\label{concatenate}
[x_0,x_1,\ldots, x_n,x_0]\times [x_0,y_1,\ldots, y_m]\longmapsto
[x_0,x_1,\ldots, x_n,x_0,y_1,\ldots, y_m]
\end{equation}
This action is continuous for the same reason presented in the proof of Corollary \ref{top}. The quotient by this action is a copy of $X$, and the projection $E(X)\lrar X$
is identified with the evaluation at the endpoint of the path; the reason being that two paths starting at $x_0$, ending at $y$ and traversed in opposite directions, give a trivial loop at $x_0$.
\end{proof}

We finally reach the main result of this paper (Theorem \ref{main}).

\bco\label{finally} There are homotopy equivalences
$\omega(X)\simeq\Omega (X)\simeq\Omega_{pl}^w(X)$.\eco

\begin{proof}
We've shown in Theorem \ref{pltos} that $\omega (X)\cong S(X)$
and $P^{thin}(X)\cong E(X)$ in the weak topology.
On the other hand, we've shown in this section that $E(X)\rightarrow X$ is a (principal) bundle with contractible total space and fiber $S(X)$.
The Puppe sequence $\Omega (X)\rightarrow S(X)\rightarrow E(X)\rightarrow X$
gives that $\Omega (X)$ is weakly equivalent to $S(X)$. By a result of Milnor,
$\Omega (X)$ is of the homotopy type of a CW complex. We need show the same is true
for $\omega (X)$ in order to see that both spaces are homotopy equivalent.
Note that each $S_n$ is of the homotopy
type of a CW complex (see Proof of Proposition \ref{contractible}). Since
$S_n(X)$ is in hCW, so is the direct limit $S(X)\cong\omega (X)$. This
establishes the equivalence $S(X)\simeq \Omega (X)$.

The homotopy equivalence $\Omega_{pl}^wX\simeq\Omega_{pl}X$ is obtained similarly by comparing path-loops fibrations. The evaluation of paths
$P_{pl}^w(X)\rightarrow X$ is a fibration in the weak topology, with contractible total space. Moreover, the identity
$P_{pl}^w(X)\rightarrow P_{pl}(X)$ is continuous and a map of path-loop fibrations. By comparison, $\Omega_{pl}^w(X)$ is weakly $\Omega_{pl}(X)$, and since these spaces have again the homotopy type of a CW, they are homotopy equivalent. Finally, $\Omega_{pl}(X)\simeq\Omega (X)$ by Lemma \ref{pspl}.
\end{proof}

\bre We only check here that $\omega (X)$ is of the homotopy type of a CW complex, and we strongly suspect it is CW on the nose.
\ere

Finally we complete the picture by discussing $LX$ the space of free loops on $X$ with the compact-open topology.
The corresponding space of piecewise loops $L_{pl}X$
comes with a filtration as in the case of $\Omega_{pl}(X)$, and we define
$\ell (X)$ to be the space of thin free loops equipped with the weak topology.

\bpr\label{freeloop} There is a homotopy equivalence $\ell (X)\simeq L(X)$.
\epr

\begin{proof} We will write $L_{pl}^w(X)$ the free PL loops on $X$ with the weak topology induced from the filtration
by vertices. Both projections $L_{pl}(X)\rightarrow X$ and
$L_{pl}^w(X)\rightarrow X$ are quasibrations by Proposition \ref{freeloopquasi}. On the other hand, $\ell (X)$ maps onto $X$ via the evaluation as well and this map is now a bundle (same proof as Lemma \ref{principal}). We have a diagram of quasifibrations (the vertical maps)
$$\xymatrix{
\omega (X)\ar[d]&\Omega^w_{pl} (X)\ar[r]\ar[d]\ar[l]& \Omega_{pl} (X)\ar[d]\\
\ell (X)\ar[d]&L^w_{pl}(X)\ar[d]\ar[r]\ar[l]& L_{pl}(X)\ar[d]\ar[u]\\
X&X\ar[r]^=\ar[l]_=&X
}$$
All spaces are in hCW (see Lemma \ref{pspl} and the proof of Proposition \ref{contractible}). Since left and right spaces are homotopy equivalent, the middle spaces are homotopy equivalent as well. But then $L_{pl}(X)\simeq L(X)$, again by Proposition \ref{freeloopquasi}.
\end{proof}


\section{Higher Thin Loop Spaces}\label{higher}

In this final section we treat the case of maps of higher dimensional spheres (cf. \cite{picken}). Here $X$ is a simplicial complex as before, $S^n$ the unit sphere in $\bbr^{n+1}$. We say $f,g: S^n\lrar X$ are thin homotopic if there is a sequence of maps
$f_1=f,\cdots , f_n=g$, such that for each $i\geq 1$, $f_i$ and $f_{i+1}$
are homotopic and the homotopy between them
is contained in their images.
We can then define the $n$-th thin loop space to be the quotient
$$\pi_n^1(X) = \{S^n\lrar X, \infty\mapsto x_0\}/_{\smaller\hbox{thin homotopy}}$$
As before there is a thin loop product which is associative with inverses and an identity (the constant loop) so that $\pi_n^1 (X)$ is a group. As in the case $n=1$, we restrict ourselves to the space of piecewise smooth or piecewise linear loops defined analogously as in Definition \ref{defpl}. In that case however the analog of Theorem \ref{homeo} fails.

The first observation here is that for $n\geq 2$, $\pi_n^1(X)$ is an \textit{abelian} group for either PL or PS maps. Indeed the homotopy between $fg$ and  $gf$ for loops in dimension at least two, as often described in a sequence of diagrams in every school book on homotopy theory, is a thin homotopy.
Secondly, a connected abelian topological group is a GEM, i.e. a space of the weak homotopy type of a product of Eilenberg-MacLane spaces (\cite{hatcher}, Theorem 4.K.6). From this we can deduce the following proposition.

\bpr Suppose all maps are piecewise smooth. Then
$\pi_2^1(S^3)$ is not of the weak homotopy type of $\Omega^2(S^3)$.
\epr

\begin{proof}
We wish to show that $\Omega^2(S^3)$ cannot be weakly a GEM $Y$. We shall recall that a weak homotopy equivalence induces an isomorphism of the singular homology groups of the spaces involved. 
We therefore need show that $\Omega^2S^3$ doesn't have the homology of the GEM below
\begin{eqnarray*}
\prod K(\pi_i(\Omega^2S^3),i)&\simeq& \prod K(\pi_{i+2}(S^3),i)\\
&\simeq& K(\bbz, 1)\times K(\bbz_2,2)\times K(\bbz_2,3)\times\prod K(G_i,i)\ \ ,\ \ i\geq 3
\end{eqnarray*}
With mod-$2$ coefficients, $H_*(\Omega^2 (S^3))$ is a polynomial algebra on a degree one generator $\iota$ (this is a result of Kudo and Araki. For a full treatment, see (\cite{browder}, Theorem 3). The mod $2$ betti numbers for the double loop space are $b_n=1$ for all $n$. For the righthand product,
$b_3=2$ generated by the fundamental class of
$K(\bbz_2,3)$ and by the product of the fundamental classes in
$K(\bbz,1)\times K(\bbz_2,2)$. This is a contradiction.
\end{proof}


\section{Appendix: Group Models for Loop Spaces}\label{models}

The loop space $\Omega (X)$ is an $H$-space which is associative only up to homotopy, and more generally an $A_\infty$ space (\cite{stasheff}). It is of interest in both homotopy theory and geometry to replace $\Omega (X)$ by a weakly homotopic space that is a topological group.
We list below the various known ways to do this.

 \subsection{Combinatorial Models}\label{milnormodel}
 Milnor in \cite{milnor} is perhaps the earliest to have exhibited such a group replacement
 for $\Omega (X)$  when $X$ is of the homotopy type of a connected countable complex. For $X$ a countable simplicial complex, consider
$S_k$ as in \eqref{sk} and define
$$M_n := \bigcup_{1\leq k\leq n}S_k/_\sim$$
where $
(x_0,\ldots, x_{i-1}, x_i, \ldots, x_n)\sim (x_0,\ldots, x_{i-1},\hat{x}_i,x_{i+1},\ldots, x_n)
$
whenever  the triple $x_{i-1}=x_i$ or $x_{i-1}=x_{i+1}$. Then
$M(X)=\lim M_n$ is the group of a universal bundle, thus it is weakly homotopy equivalent to $\Omega (X)$. Our model in this paper $S(X)$ has the advantage of being homeomorphic to the thin loop space for piecewise linear loops. Note that there is a homomorphism of topological groups
$M(X)\rightarrow\omega (X)$.

Another model appears in \cite{bc} in the case $X$ is a Riemannian manifold. The loop space model by Bahri and Cohen is obtained by considering composable small geodesics.

 \subsection{Simplicial (Kan) Model}. Let $X$ be a based space which is of the homotopy type of a CW complex. Take the simplicial total singular complex $SX$ of $X$ which is a based simplical set. The Kan loop group $GSX$ of $SX$ is then a simplicial group and its geometric realization $|GSX|$ is a topological group object in the category of compactly generated spaces. This is a model for the loop space. This is discussed in \cite{klein}.

 The generality of this model allows us to deduce for example that
 any strictly associative monoid (or $H$-space) such that $\pi_0$ is a group can be ``rigidified" to a topological group; that is is weakly homotopy equivalent to a topological group. Indeed for such spaces, Dold and Lashof show the existence of a classifying space $B_X$ and a weak homotopy equivalence $X\lrar\Omega B_X$. One then uses the above construction to replace the loop space
up to equivalence by a topological group.

Note that P. May has another model that produces this time a topological monoid weakly equivalent to any $A_\infty$-space (Theorem 13.5 of \cite{may}).

\subsection{The ``Geometric Model"}
Lefshetz (\cite{lefshetz}, Chap. V, \S 4) and Kobayashi introduce the earliest known thin loop models for piecewise smooth loops on smooth manifolds. Teleman in \cite{teleman} introduces
a similar model which is a quotient of $P(M,x_0)$, the space of piecewise smooth paths starting at $x_0$, by three equivalence relations and uses this model to classify bundles.
Gajer introduces a similar model in (\cite{gajer}, \S1.3), which
with the first relation ensuring that in the quotient space $G(M)$, the product of loops is associative, $\sim_2$ ensures that the class of the constant map is the identity in $G(M)$, and the final relation $\sim_3$ implies that $[\gamma^{-1}]$ is a genuine inverse for $[\gamma]$, much as in our \S\ref{thin}. However, the claim that $P^g(X)$ is a model for $\Omega (X)$ is not fully justified (see our remark \ref{thepscase}).

\subsection{The Quillen Model} Let $X=G$ be a simply-connected Lie group with simple Lie algebra, and let $G_\bbc$ be its complexification.
An ``algebraic loop" in $G_\bbc$ is a regular map $\bbc^*\lrar G_\bbc$ (i.e. a morphism of algebraic varieties). The set of all such
loops is denoted by $\tilde{G}_\bbc$. One then sets
$$\Omega_{alg}G := \{f\in \tilde{G}_\bbc\ |\ f_{|S^1}\in\Omega G\}$$
This is an algebraic group.
By a Theorem of Quillen and Garland-Raghunathan,  the restriction map to $S^1$ yields a homotopy equivalence $\Omega_{alg}G\simeq\Omega G$ (one reference is \cite{mitchell}).

\end{document}